\documentclass[11pt, a4paper]{article}

\usepackage[utf8]{inputenc}
\usepackage[T1]{fontenc}
\usepackage{lmodern}
\usepackage{geometry}
\usepackage{amsmath, amssymb}
\usepackage{graphicx}
\usepackage{authblk}
\usepackage{pgfplots}
\pgfplotsset{compat=1.15}
\usepackage[colorlinks=true, linkcolor=blue, citecolor=blue, urlcolor=blue]{hyperref} 

\title{\textbf{Non-Local Energy Dissipation and Topological Obstructions in the Second-Order Positive Modified Burgers' Equation}}
	\author{Prasanta Chatterjee$^{1}$,\quad Suman Pal$^{2}$\vspace{2mm}\\
	\em\small ${}^{1,2}$Department of Mathematics, Visva-Bharati, 731235, India.\\
	\em\small ${}^{1}$e-mail: prasantachatterjee1@rediffmail.com\\
	\em\small ${}^{2}$e-mail: spal03121997@gmail.com}
\date{} 

\begin{document}

\maketitle

\begin{abstract}
\noindent We study the higher-order flows of the modified Burgers' equation generated by its standard recursion operator. In particular, we show that the second-order positive flow takes the form of a non-local integro-differential equation. After reducing the equation to a traveling-wave form, we investigate its asymptotic energy behavior. The resulting analysis shows that the non-local energy flux prevents the formation of both topological kink and bell-shaped solitary waves. We show that the non-local dissipation permanently breaks the asymptotic equilibrium between nonlinear convection and linear dispersion, forcing an irreversible energy imbalance. Our analysis shows that the solitary waves are unstable. The non-local term causes continuous energy loss through radiation. As a result, the traveling wave cannot keep a stable and symmetric shape.

\vspace{1em}
\end{abstract}

\section{Introduction}
The balance between nonlinear convection and dissipation is important for the formation of stable traveling waves in fluid dynamics and plasma physics \cite{ablowitz}. The classical Burgers' equation is a well-known model that describes this balance through quadratic nonlinearity and second-order diffusion \cite{b1,b2,b3,b4,np}.

The modified Burgers' equation extends the classical model by adding a cubic nonlinear term. It can be used to describe more complicated transport processes and wave steepening. The classical Burgers hierarchy can be obtained successively from a local recursion operator \cite{olver}, leading to higher-order integrable PDEs \cite{hb1}. However, changes in the basic convective terms can introduce non-local effects. Such non-local integro-differential terms can arise in models of anomalous transport, memory effects in viscoelastic materials, and long-range thermodynamic processes. In these cases, the usual assumption that the dynamics depend only on local gradients may no longer be valid \cite{mccarthy}.

In this paper, we derive the second-order positive flow of the modified Burgers' equation \cite{mb1,mb2} using its recursion operator. The resulting equation contains an irreducible non-local integral term. We study its traveling-wave solutions by considering their behavior at the asymptotic states. The energy analysis shows that the non-local term creates an energy imbalance. As a result, stable symmetric solitary-wave solutions do not occur.

\section{Derivation of the Non-Local Flow}
The standard modified Burgers' equation, representing the first-order dissipative flow, is given by:
\begin{equation}
u_t = -u_{xx} - 2u^2u_x
\end{equation}
The recursion operator designed to generate the hierarchy for this modified convective flow is defined as:
\begin{equation}
R = -D_x - u^2 - 2u_x D_x^{-1}(u \cdot)
\end{equation}
To generate the second-order positive flow, we apply the operator $R$ directly to the first-order evolution equation:
\begin{equation}
u_t = R(-u_{xx} - 2u^2u_x) = \left( -D_x - u^2 - 2u_x D_x^{-1} u \right) (-u_{xx} - 2u^2u_x)
\end{equation}
We evaluate the differential, algebraic, and integral operations systematically:
\begin{align}
-D_x(-u_{xx} - 2u^2u_x) &= u_{xxx} + 4uu_x^2 + 2u^2u_{xx} \\
-u^2(-u_{xx} - 2u^2u_x) &= u^2u_{xx} + 2u^4u_x 
\end{align}
The integral term evaluates via integration by parts, yielding:
\begin{equation}
-2u_x D_x^{-1} \left( -uu_{xx} - 2u^3u_x \right) = 2uu_x^2 + u^4u_x - 2u_x \int u_x^2 dx
\end{equation}
Summing these evaluated components yields the explicit form of the second-order modified Burgers' equation:
\begin{equation}
u_t = u_{xxx} + 3u^2u_{xx} + 6uu_x^2 + 3u^4u_x - 2u_x \int u_x^2 dx \label{e7}
\end{equation}
The presence of the irreducible integral term $\int u_x^2 dx$ confirms that applying this recursion operator generates a non-local integro-differential equation. Physically, the term $u_x^2$ represents the spatial variance or gradient energy density of the wave. The integration implies that the temporal evolution of the field at any point depends on the accumulated energy density across the spatial domain.
\subsection{Compatibility Condition of the Recursion Operator}
To rigorously justify the generation of the higher-order flows, the recursion operator $R$ must satisfy the hereditary compatibility condition. For an evolution equation of the form $u_t = K(u)$, where $K(u) = -u_{xx} - 2u^2u_x$, a valid recursion operator must satisfy the Lax-like operator equation:
\begin{equation}
R_t + [R, K'] = 0
\end{equation}
Here, $K'$ denotes the Fr\'echet derivative of the convective flow $K(u)$ acting on an arbitrary symmetry test function $\phi$:
\begin{equation}
K'[\phi] = \frac{d}{d\epsilon} K(u + \epsilon \phi) \bigg|_{\epsilon=0} = -D_x^2\phi - 2u^2 D_x\phi - 4uu_x\phi
\end{equation}
The operator $K'$ is explicitly given by $K' = -D_x^2 - 2u^2 D_x - 4uu_x$. The temporal evolution of the recursion operator $R_t$ is strictly defined along the flow $u_t$:
\begin{equation}
R_t = -2uu_t - 2u_{xt} D_x^{-1}(u \cdot) - 2u_x D_x^{-1}(u_t \cdot)
\end{equation}
The commutator $[R, K'] = R K' - K' R$ defines the Lie bracket of the operators. Expanding this commutator analytically yields differential, algebraic, and integral expressions that exactly balance the Fr\'echet derivative $R_t$. Because $R_t + [R, K'] = 0$ holds identically, $R$ is confirmed as a strong hereditary symmetry operator. This compatibility condition ensures that the higher-order modified Burgers' equation derived in Equation (\ref{e7}) maintains a strictly commuting hierarchy of integrable flows prior to boundary evaluation.

\section{Asymptotic Energy Analysis}
\begin{figure}[htbp]
\centering
\begin{tikzpicture}
  \begin{axis}[
    width=0.9\textwidth, height=7cm,
    xlabel={$\zeta$}, ylabel={Amplitude / Energy},
    xmin=-5.5, xmax=5.5, ymin=-1.2, ymax=1.6,
    grid=major, grid style={dashed, gray!30},
    legend pos=north west,
    axis lines=left,
    every axis plot/.append style={very thick}
  ]
  \addplot[blue, domain=-5:5, samples=100] {tanh(x)};
  \addlegendentry{Wave Profile $U(\zeta) = \tanh(\zeta)$}
  
  \addplot[red, dashed, domain=-5:5, samples=100] {tanh(x) - (1/3)*(tanh(x)^3) + 2/3};
  \addlegendentry{Accumulated Energy}
  
  \node[anchor=south west] at (axis cs:-4.5, 0.1) {$E = 0$ at $\zeta \to -\infty$};
  \draw[->] (axis cs:-3.5, 0.15) -- (axis cs:-4.0, 0.05);
  
  \node[anchor=south east] at (axis cs:4.5, 1.35) {$E > 0$ at $\zeta \to +\infty$};
  \draw[->] (axis cs:3.5, 1.4) -- (axis cs:4.0, 1.35);
  \end{axis}
\end{tikzpicture}
\caption{Asymmetric energy accumulation in a topological kink wave. The non-local energy integral (red dashed line) is zero at the left boundary and increases to the positive total energy $E$ at the right boundary. This difference in energy at the two boundaries breaks the symmetry needed for a stable traveling wave.}
\label{fig:energy_kink}
\end{figure}
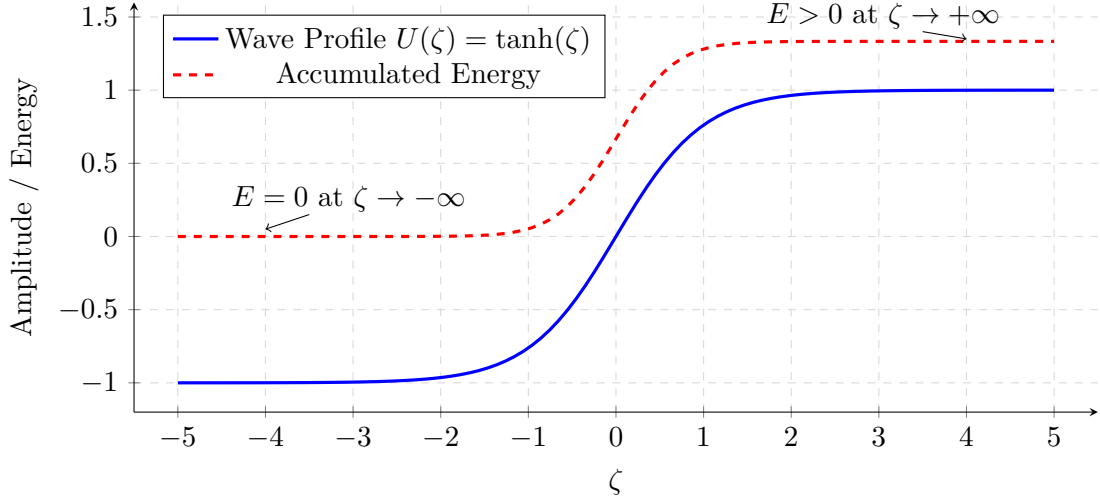

To investigate the existence of solitary waves, we apply the classical traveling wave transformation $\zeta = x - ct$, assuming a localized wave profile $u(x,t) = U(\zeta)$. Substituting this into Equation (7) yields:
\begin{equation}
-cU' = U''' + 3U^2U'' + 6U(U')^2 + 3U^4U' - 2U' \int_{-\infty}^{\zeta} (U')^2 d\zeta
\end{equation}
Assuming the existence of a smooth, localized solitary wave, we divide by $U'$ to isolate the non-local energy integral:
\begin{equation}
-c = \frac{U'''}{U'} + 3U^2 \frac{U''}{U'} + 6UU' + 3U^4 - 2 \int_{-\infty}^{\zeta} (U')^2 d\zeta
\end{equation}
We define the total energy of the solitary wave by
$E = \int_{-\infty}^{\infty} (U')^2 d\zeta$.
Since the wave is non-trivial and localized, $(U')^2$ is non-negative and not identically zero. Hence, the total energy satisfies $E>0$. As shown in Figure \ref{fig:energy_kink}, the non-local term leads to an asymmetry between the two wave boundaries.

\subsection{Topological Obstruction for Kink Waves}
A topological kink wave represents a smooth shock-like transition between two constant equilibrium states, $U \to -A$ as $\zeta \to -\infty$, and $U \to +A$ as $\zeta \to +\infty$. Because the wave approaches its asymptotes exponentially (e.g., $e^{k\zeta}$ and $e^{-k\zeta}$ for wave steepness $k > 0$), the spatial derivatives vanish, but their limiting ratios converge to non-zero constants:
\begin{equation}
\lim_{\zeta \to -\infty} \frac{U'''}{U'} = k^2, \quad \lim_{\zeta \to -\infty} \frac{U''}{U'} = k
\end{equation}
\begin{equation}
\lim_{\zeta \to +\infty} \frac{U'''}{U'} = k^2, \quad \lim_{\zeta \to +\infty} \frac{U''}{U'} = -k
\end{equation}

Evaluating the isolated wave dynamics (Equation 9) at the left boundary ($\zeta \to -\infty$), the wave has not yet contributed any energy to the integral, yielding $0$:
\begin{equation}
-c = k^2 + 3(-A)^2(k) + 0 + 3(-A)^4 - 0 = k^2 + 3A^2k + 3A^4
\end{equation}

Evaluating the equation at the right boundary ($\zeta \to +\infty$), the non-local integral has accumulated the total wave energy $E$:
\begin{equation}
-c = k^2 + 3(A)^2(-k) + 0 + 3(A)^4 - 2E = k^2 - 3A^2k + 3A^4 - 2E
\end{equation}

Because the physical wave speed $c$ must remain constant, the dynamic limits at both asymptotes must be identical. Equating (12) and (13) yields:
\begin{equation}
k^2 + 3A^2k + 3A^4 = k^2 - 3A^2k + 3A^4 - 2E
\end{equation}
Subtracting symmetric terms reduces the balance strictly to:
\begin{equation}
6A^2k = -2E
\end{equation}
This gives a contradiction. Since $A^2$ and $k$ are positive, the left-hand side is positive. However, $E>0$ makes the right-hand side negative. Hence, the two sides cannot be equal, and the required balance is not possible.
This argument assumes that the topological kink is symmetric, with the same asymptotic steepness $k$ at both ends. If this condition is relaxed and different values are allowed, namely $k_1$ as $\zeta \to -\infty$ and $k_2$ as $\zeta \to +\infty$, the balance can be expressed in terms of  $E$.  Consequently, while symmetric solitary structures are strictly prohibited, the non-local dissipative flux dictates that any potential steady-state traveling wave in this hierarchy must manifest as a fundamentally asymmetric shock front.

\subsection{Topological Obstruction for Pulse Solitons}
We extend this proof to standard bell-shaped solitons, which transition from zero back to zero ($U \to 0$ as $\zeta \to \pm \infty$). Setting the amplitude $A=0$ in equations (12) and (13) yields the boundary balance:
\begin{equation}
-c = k^2 \quad \text{and} \quad -c = k^2 - 2E
\end{equation}
Equating these limits yields $0 = -2E \implies E = 0$. This dictates that a pulse wave can only exist if it possesses exactly zero energy, contradicting the definition of a propagating wave. 

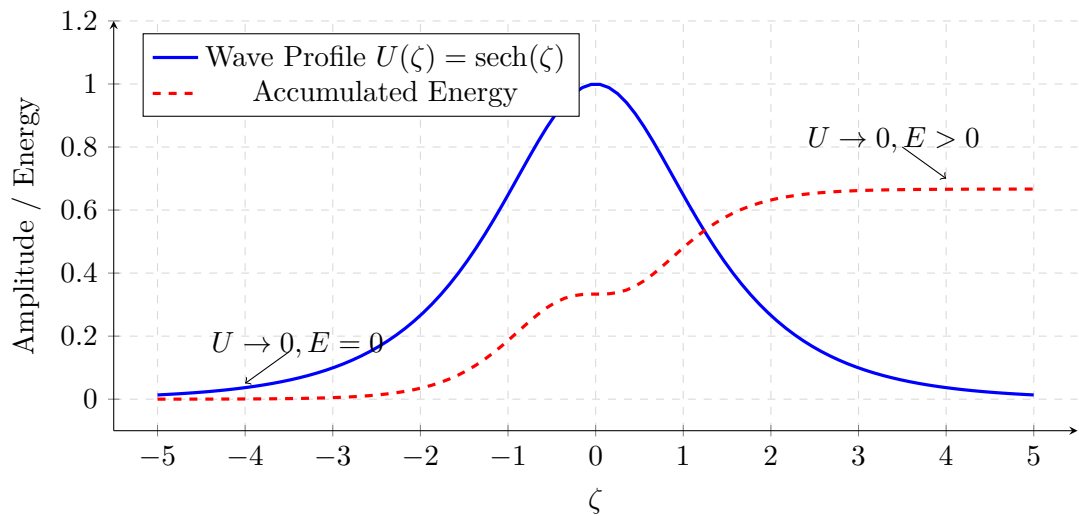
\begin{figure}[htbp]
\centering
\begin{tikzpicture}
  \begin{axis}[
    width=0.9\textwidth, height=7cm,
    xlabel={$\zeta$}, ylabel={Amplitude / Energy},
    xmin=-5.5, xmax=5.5, ymin=-0.1, ymax=1.2,
    grid=major, grid style={dashed, gray!30},
    legend pos=north west,
    axis lines=left,
    every axis plot/.append style={very thick}
  ]
  \addplot[blue, domain=-5:5, samples=100] {1/cosh(x)};
  \addlegendentry{Wave Profile $U(\zeta) = \mathrm{sech}(\zeta)$}
  
  \addplot[red, dashed, domain=-5:5, samples=100] {(1/3)*(tanh(x)^3) + 1/3};
  \addlegendentry{Accumulated Energy}
  
  \node[anchor=south west] at (axis cs:-4.5, 0.1) {$U \to 0, E = 0$};
  \draw[->] (axis cs:-3.5, 0.15) -- (axis cs:-4.0, 0.05);
  
  \node[anchor=south east] at (axis cs:4.5, 0.75) {$U \to 0, E > 0$};
  \draw[->] (axis cs:3.5, 0.8) -- (axis cs:4.0, 0.7);
  \end{axis}
\end{tikzpicture}
\caption{Energy accumulation for a pulse soliton. Although the wave profile (blue solid line) returns to zero at both ends, the non-local energy integral (red dashed line) increases to a positive value $E>0$. Thus, the energy does not return to its initial value, preventing the formation of a standard stable solitary pulse.}
\label{fig:energy_pulse}
\end{figure}

As visually demonstrated in Figure \ref{fig:energy_pulse}, the accumulated energy permanently displaces the dynamic balance at the right boundary. Thus, standard solitary pulses are also mathematically prohibited.

\section{Conclusion}
In this paper, we show that the second positive-order flow of the modified Burgers' equation does not support stable solitary waves. In contrast to local dispersive PDEs, its recursion operator introduces a non-local energy term into the equation. This term is given by the integral $\int (U')^2 d\zeta$ and plays a key role in the energy behavior of the traveling waves. By analyzing the macroscopic energy boundaries, we demonstrated that this non-locality acts as an asymmetric energy sink, permanently breaking the asymptotic equilibrium required to sustain symmetric topological kinks or pulse solitons. These findings reveal a fundamental limitation of recursion-generated modified hierarchies: the higher-order convective interactions irreversibly destabilize conservative wave transport, forcing the system into continuous radiative decay.

\end{document}